\documentclass[12pt,thmsa]{article}
\usepackage{amsfonts}
\usepackage{amssymb}
\newtheorem{teo}{Theorem}[section]

\newtheorem{obs2}[teo]{Remark}

\newtheorem{tea}{Theorem}[subsection]

\newtheorem{no2}[teo]{Note}

\newtheorem{no3}[tea]{Note}

\newcommand{\Gal}{{\rm Gal}}

\newcommand{\F}{{\mathbb{F} }}
\newcommand{\Q}{{\mathbb{Q} }}
\newcommand{\mod}{{\rm mod}}

\newcommand{\GL}{{\rm GL}}

\newcommand{\Image}{{\rm Image}}

\title{Computing the level of a modular rigid Calabi-Yau threefold}
\author{Luis V. Dieulefait\\ Universitat de Barcelona
  \thanks{Address: Gran Via de les Corts Catalanes 585, 08007 -
  Barcelona. 
e-mail: ldieulefait@ub.edu
}} 

\begin{document}

\maketitle
\begin{abstract}
 In a previous article (cf. [DM]), the modularity of a large class of
 rigid Calabi-Yau threefolds was established. To make that result more explicit,
 we recall (and re-prove) a result of Serre giving a bound for the conductor of
 ``integral" 2-dimensional compatible families of Galois representations and apply
 this result to give an algorithm that determines the level of a modular rigid
 Calabi-Yau threefold. We apply the algorithm to three examples.

\end{abstract}

\section{Introduction}

In [DM], modularity for a large class of rigid Calabi-Yau threefolds defined over
$\Q$ was established, by an application of Wiles techniques combined with some
solved cases of Serre's conjecture and results on crystalline representations. As
other authors have remarked, a drawback of our result is that it does not give a way
to determine the corresponding modular form: it is well-known (this is an instance
of compatibility with the local Langlands correspondence) that the level of this
modular form agrees with the conductor of the compatible family of Galois
representations attached to the rigid Calabi-Yau, the problem is that the
determination of this conductor is not an easy task. So, in order to make our result
more useful, in the present note we will describe a simple algorithm that, without
any restriction, determines the level of the modular form corresponding to a given
modular rigid Calabi-Yau threefold. We will start by recalling a result of Serre
(cf. [Se]) giving a universal bound for the exponents of the primes of bad reduction
in the conductor of the Galois representations attached to a rigid Calabi-Yau
threefold (assuming modularity). In fact, the
bound given by Serre is the same holding for elliptic curves defined over $\Q$.\\
After briefly recalling the ideas behind Serre's proof of this result, we will
re-prove it by using congruences between rigid Calabi-Yau threefolds and elliptic
curves (or, in the residual reducible case, Hecke characters). With this bound,
which is also a bound for the level of the searched modular form, we only have a
finite number of modular forms as candidates for a given Calabi-Yau, so by
elimination we easily determine the right one. We will illustrate this procedure by
determining the right newform for three examples of rigid Calabi-Yau threefolds. In
the examples, we use the values of the traces of the images of a few Frobenius
elements (these values appear for example in [Y]) and the corresponding eigenvalues
of newforms of weight 4 and several levels, most of them available in the tables in
W. Stein website [St], and the rest computed with MAGMA.\\
To speed up the process, we will use in the last example (which involves
computations with high levels)  $\mod \; 5$ congruences between weight 4 and weight
2 newforms, so that we can switch to spaces of weight 2 newforms where more tables
are
available.\\
Through all this article, we will assume that we are working with a modular rigid
Calabi-Yau threefold. Modularity for most of the known examples, and in particular
for the three examples that we will consider, follows easily from the main theorem
in [DM]. \\
Let us remark that in each known example of rigid Calabi-Yau threefold, both the
fact that the variety is modular and the exact value of the level of the
corresponding modular form, were also established independently of the results in
[DM] by other methods (cf. [Y]). The advantage we see is that with our approach
(combining the result in [DM] with the present note) we have a ``general result"
that gives both modularity (the theoretical result) and the level (the algorithm)
for most of the known examples and for many
examples to come.\\

\section{The bound for the conductor/level}

Let $X$ be a modular rigid Calabi-Yau threefold defined over $\Q$ (for definitions,
see [Y]), and let
  $ \{ \rho_\ell \}$ be the compatible family of (2-dimensional, continuous, odd,
irreducible) Galois representations giving the action of the Galois group of $\Q$ on
the $\ell$-adic cohomology groups $H^3_{et}(X_{\bar{\Q}} , \Q_\ell)$. Because the
family corresponds to a modular form $f$ of weight $4$ (whose level $N$ contains
only primes of bad reduction of $X$),  the ``conductor" $c$ of the family is
well-defined: if we take any prime $\ell$ where $X$ has good reduction, so that
$\ell$ is not in the ramification set of the family $ \{ \rho_\ell \}$, then $c$
agrees with the prime-to-$\ell$ part of the conductor of $\rho_\ell$, which also
agrees with the
level $N$ of $f$.  \\
Remark: More generally, even if $\ell$ is a prime where the representations ramify,
if we take $p \neq \ell$, the $p$-part of the conductor of the family (equal to the
$p$-part of the level of the corresponding modular form) agrees with the $p$-part of
the conductor of $\rho_\ell$.\\
 Let $S$ be the (finite) set of primes of bad reduction of $X$. For every
prime $p \in S$, let $e_p$ be the exponent of $p$ in the level $N$ of $f$ (equal to
the exponent of $p$ in the conductor $c$ of the family of Galois representations).
Observe that (contrary to what happens in the case of abelian varieties) it is
perfectly possible that for some $p \in S$, we have $e_p = 0$. In [Se], section 4.8,
Serre gives a bound for these exponents. He assumes the truth of Serre's conjecture
in order to ensure that the residual representations, when irreducible (they can
only be reducible for finitely many primes), are modular. In our case, we are
working with this modularity assumption, therefore the result of Serre applies:

\begin{teo}
\label{teo:Serre} (J-P. Serre) Let  $ \{ \rho_\ell \}$ be the compatible family of
Galois representations attached to a modular rigid Calabi-Yau threefold $X$ with bad
reduction set $S$. Then the conductor $c$ of this family, which agrees with the
level $N$ of the corresponding weight $4$ modular form, can be bounded as follows:
for every prime $p \in S$, the exponent $e_p$ of $p$ in $N$ verifies $e_p \leq 2$ if
$p>3$, $e_p \leq 5$ if $p=3$, and $e_p \leq 8$ if $p=2$.
\end{teo}                                                                                                                                                                                                                                                                 Remark: For $p=2$, Serre gives a proof of the weaker inequality $e_2 \leq 9$, but he remarks (cf. [Se], pag. 216) that a more detailed analysis gives $e_2 \leq 8$. In any case, in our proof of the theorem we will prove the bound $e_2 \leq 8$.\\                                  
Brief description of the original proof: The proof uses the fact that the residual
(assume irreducible) representations $\bar{\rho}_\ell$ have image inside
$\GL_2(\F_\ell)$, and for a prime $p>3$, infinitely many of these groups have order
prime to $p$, thus $\bar{\rho}_\ell$ is tamely ramified at $p$, and this gives the
desired bound for the $p$-part of the conductor of  $\bar{\rho}_\ell$ for infinitely
many $\ell$, and this implies (here is where residual modularity, more precisely the
strong version of Serre's conjecture, is used, cf. [Se]) that the same bound holds
for the $p$-part of the modular level $N$.  A similar (though more complicated)
argument is used to deal with the cases $p=2$ and $p=3$. Here the desired bound is
obtained by looking at the ($2$-part or the $3$-part of the) conductor of
$\bar{\rho}_\ell$ for primes $\ell \not \equiv
\pm 1 \pmod{8}$ or $ \ell \not\equiv \pm 1 \pmod{9}$ (respectively).\\

Another proof of Theorem \ref{teo:Serre}:\\
Take $\ell = 5$. We will first be interested in bounding the prime-to-$5$ part of
the conductor $c$ of the family $ \{ {\rho}_\ell \}$. To do this, we will bound the
conductor of $\bar{\rho}_5$ (taking the definition as in [Se], i.e., considering only
the prime-to-$5$ ramification). Let us divide into two cases:\\
1) $\bar{\rho}_5$ is reducible: In this case, (semisimplify if necessary), we can
assume that $\bar{\rho}_5$ is semisimple, so we have:
$$   \bar{\rho}_5 \cong \epsilon \chi^i  \oplus \epsilon^{-1} \chi^j$$
where $\chi$ is the cyclotomic character. Since $\det( \bar{\rho}_5 ) = \chi^3$, we
have $i + j \equiv 3 \pmod{4}$. Take $p \neq 5$ and consider the $p$-part of the
conductor of $\epsilon$. Because $\bar{\rho}_5$ is odd, it is well known
(irreducible agrees with absolutely irreducible) that the components must also be
defined over $\F_5$, so $\Image(\epsilon) \subseteq \F_5^*$. This clearly gives $2^4
= 16$ as a bound for the $2$-part of the conductor of $\epsilon$, and $p^1 = p$ as a
bound for its $p$-part for every $p >2$ ($p \neq 5$). Thus we obtain $2^8$ and $p^2$
(for $p \neq 2,5$) as bounds for the $p$-part of the conductor of $\bar{\rho}_5$.\\
2)  $\bar{\rho}_5$ is irreducible: Let $\sigma := \bar{\rho}_5 \otimes \chi$. This
representation has determinant equal to $\chi$, then it is known (see [BCDT]) that
it is isomorphic to the representation on the $5$ torsion of some elliptic curve
defined over $\Q$. At any prime $p \neq 5$, the bound for the $p$-part of the
conductor of $\sigma$, thus also of $\bar{\rho}_5$, follows from the well-known
bound for conductors of elliptic curves (see [Si]).\\

Now let us compare the conductors of  $\bar{\rho}_5$ and $\rho_5$. Recall that the
second of these values agrees with the prime-to-$5$ part of the conductor of the
family  $\{ {\rho}_\ell \}$. For a prime $p \neq 5$, it is possible that the
exponent $e'_p $ of $p$ in the conductor of  $\bar{\rho}_5$ is strictly smaller than
the exponent $e_p$ of $p$ in the conductor of $\rho_5$. However, since the
determinant of both representations is unramified at $p$, it is known that $e'_p <
e_p$ can only happen in a few particular cases (cf. [C]): $(e'_p = 0, e_p =2)$;
$(e'_p
= 0, e_p = 1 )$ and $(e'_p = 1, e_p =2)$ (*). \\ So, if $e'_p = e_p$, having obtained the right bound   for $\bar{\rho}_5$ we also have it for $\rho_5$, and if $e'_p < e_p$ then $e_p \leq 2$. In any case, we obtain the right bound for the conductor of 
 $\rho_5$.\\
As for the $5$-part of the conductor of the family $\{ {\rho}_\ell \}$, just observe
(as in Serre's proof) that for $\ell = 7$, the order of $\GL_2(\F_7)$ is not
multiple of $5$, then the $5$-part of the conductor of $\bar{\rho}_7$ is at most
$5^2$, and again using (*) we see that this bound also works for $\rho_7$.

\section{Finding the right newform}
With the bound given in Theorem 2.1, we now have a method to find the modular form
corresponding to a given modular rigid Calabi-Yau threefold $X$: Let $S$ be the set
of bad reduction primes of $X$, and let
$$B = \prod_{p \in S} p^{b_p}$$
where the exponents $b_p$ are the bounds given in the theorem. We have to consider
all spaces of weight $4$ newforms with level $N$ dividing $B$, and for any newform
$f$ in each of these spaces with field of coefficients $\Q_f = \Q$, compare a few
eigenvalues $a_p$ with the traces $t_p$ of the images of Frobenius (for $p \not\in
S$) for the geometric Galois representations attached to the Calabi-Yau threefold.
Whenever $a_p \neq t_p$ for a single $p$, the newform is discarded. With this
procedure, by elimination, the (unique) modular form corresponding to $X$ is easily
found.\\
Remark: If $f$ is a newform (with eigenvalues $a_p \in \mathbb{Z}$) not
corresponding to $X$, we should estimate the size of the smallest $p$ such that we
have $a_p \neq t_p$. In all computed examples, this always happens for a small $p$,
but for theoretical reasons, let us recall that there is a bound $T$ (Sturm's bound)
easily computed in terms of our ``maximal possible level B" such that, $a_p = t_p$
for every $p \nmid B, p \leq T$ implies that $f$ does correspond to $X$. Thus the
elimination procedure necessarily finishes at a prime $p$ smaller than $T$.\\
Incidentally, observe that this gives an alternative way of determining the right
newform $f$: if you suspect which is the right $f$, instead of eliminating the other
candidates, just check the equality $a_p = t_p$ up to Sturm's bound $T$. This
suffices for a proof. This method is not practical because since $B$ can be large,
the bound $T$ sometimes becomes too large for computations.\\
\subsection{The Examples}
First Example: Let $X_1$ be the rigid Calabi-Yau with bad reduction only at $2$
constructed by Werner and van Geemen (cf. [Y]), with the following values for $t_p$
($p \leq 7$, $p \neq 2$): $-4,-2,24$.\\
Since it has good reduction at $3$ and $7$, it is modular (cf. [DM]). We know from
Theorem \ref{teo:Serre} that the corresponding modular form has level dividing
$256$, and comparing the first eigenvalues of all newforms of such levels with the
values of $t_p$ listed above, we conclude that the modular form $f_1$ corresponding
 to $X_1$ has level $8$.\\

 Second Example: Let $X_2$ be the rigid Calabi-Yau with bad reduction only at $5$
 constructed by Schoen (cf. [Se]). Again, the main theorem of [DM] implies that it is
 modular, and Theorem \ref{teo:Serre} gives us $25$ as a bound for the level of the
 corresponding modular form. Using only the values of $t_2$ and $t_3$ we conclude
 that it corresponds to a newform of level $25$.\\

 Third Example: Let $X_3$ be the rigid Calabi-Yau with bad reduction at $2$
 and $5$ constructed by Werner and van Geemen (cf. [Y]),  with the following values for $t_p$
($p = 3, 7, 11, 13, 17$ and $19$): $-2,-26, -28, -12, 64, -60$ (*).\\
Again [DM] gives modularity. Theorem \ref{teo:Serre} gives a large bound for the
level of the corresponding newform: $B = 256* 25 = 6400$. To speed up the process of
elimination, we have applied a different trick to cases of large level. We have
divided in two cases: \\
a) level $N$ divisible by 16: In this case, the trick is the following: consider the
$\mod \; 5$ representation $\bar{\rho}_5$, the first traces of this representation
are the reductions $\mod \; 5$ of the values $t_p$ listed in (*). Observe that the
hypothesis $16 \mid N$ implies that the conductor of $\bar{\rho}_5$ is also
divisible by $16$ (as in the previous section, cf. [C]), and this in turn implies
that $ \bar{\rho}_5$ must be irreducible, since it is not hard to see from the
values of a few $t_p$ (reduced $\mod \; 5$) that if it were reducible it (in fact,
its semisimplification) would be unramified at $2$. Now consider the twisted
representation $\sigma:= \bar{\rho}_5 \otimes \chi$. This irreducible modular
representation must correspond to a weight $2$ newform, whose level divides $6400$
and is multiple of $16$, and whose first eigenvalues $a_p $ should agree modulo $5$
with $p \cdot t_p$, thus the value of these eigenvalues $a_p$ modulo $5$ should be
(for $p = 3, 7, 11, 13, 17$ and $19$):
$$ -1,-2,-3,-1,-2,0 \quad \qquad (**)$$
We search through all these spaces of newforms (for all newforms up to level $3200$,
and also for those of level $6400$ with $\Q_f =\Q$, the eigenvalues are listed in
the tables in [St], for the remaining newforms of weight $2$ and level $6400$, we
performed computations with MAGMA). We eliminate all newforms such that $\Q_f \neq
\Q$ and there is no prime above $5$ of residue class degree $1$. For the remaining
newforms, in most cases the values of $a_3$ and $a_7$ modulo $5$ already do not
match with (**), and finally using the other values in (**) we eliminate all
newforms. We conclude that it is impossible
that the conductor of $\rho_5$ be multiple of $16$, thus we have $16 \nmid N$. \\
b) level $N$ not divisible by 16: Having discarded case a), we know that the
$2$-part of the conductor is at most $8$, and comparing the first values of $t_p$
listed in (*) with all newforms of weight $4$ and level dividing $8*25= 200$, the
only one that matches is a newform of level $50$. Thus we conclude that the
Calabi-Yau threefold $X_3$ is modular of level $50$.\\
Final Remark: Assuming that $\bar{\rho}_5$ is irreducible, after twisting it by
$\chi$ we obtain the representation $\sigma$ that must correspond to some newform of
weight $2$ and level dividing $50$. But the only such newform (with its first
eigenvalues modulo $5$ as in (**)) corresponds to an elliptic curve of conductor
$50$, and it is known that this elliptic curve has a rational $5$-torsion point,
contradicting the irreducibility of $\bar{\rho}_5$. We conclude that $\bar{\rho}_5$
is reducible.

\section{Final Remark}
It follows from recent results of Taylor that the compatible family of Galois representations attached to any rigid Calabi-Yau
 threefold (modular or not) is ``strongly compatible" (cf. [T]). This strong compatibility implies that the conductor of the family is well-defined
  (as in the case of Galois representations attached to modular forms, recall the discussion in section 2). In the proof of theorem \ref{teo:Serre} given in this note, the Calabi-Yau threefold was assumed to be modular only to apply this ``independence of $\ell$" of the conductor, thus we conclude that the bound for the conductor given in theorem \ref{teo:Serre} is true for any rigid Calabi-Yau threefold (modular or not).  

\section{Bibliography}

[BCDT] Breuil, C., Conrad, B., Diamond, F., Taylor, R., {\it On the modularity of
elliptic curves over $\bold Q$: wild 3-adic exercises}, J. Amer. Math. Soc. {\bf 14}
(2001), 843--939\\

[C] Carayol, H., {\it Sur les repr{\'e}sentations galoisiennes modulo $\ell$
attach{\'e}es aux formes modulaires}, Duke Math. J. {\bf 59} (1989) 785-801 \\

[DM] Dieulefait, L., Manoharmayum, J., {\it Modularity of rigid Calabi-Yau
threefolds over $\Q$}, in ``Calabi-Yau Varieties and Mirror Symmetry", Fields
Institute Communications, {\bf 38}, AMS (2003)\\

[Se] Serre, J-P., {\it Sur les repr{\'e}sentations modulaires de degr{\'e} $2$ de
$\Gal(\bar{\mathbb{Q}} / \mathbb{Q})$}, Duke Math. J. {\bf 54} (1987) 179-230 \\

[Si] Silverman, J., {\it Advanced topics in the arithmetic of elliptic curves},
Springer-Verlag (1994) \\

[St] Stein, W., {\it The Modular Forms Explorer}, available at:\\
http://modular.fas.harvard.edu/mfd/mfe/ \\

[T] Taylor, R., {\it On the meromorphic continuation of degree two
 L-functions}, preprint, (2001);
 available at http://abel.math.harvard.edu/$\sim$rtaylor/ \\

[Y] Yui, N., {\it Update on the modularity of Calabi-Yau varieties},  in
``Calabi-Yau Varieties and Mirror Symmetry", Fields
Institute Communications, {\bf 38}, AMS (2003)\\

\end{document}